\documentclass[11pt,a4paper,reqno]{amsart}
\usepackage{amsmath}
\usepackage{amssymb,latexsym}
\usepackage[dvips]{graphicx}
\usepackage{color}
\usepackage{cite}
\usepackage{amsthm}

\newtheorem{theorem}{Theorem} [section]
\newtheorem{lemma}{Lemma} [section]
\newtheorem{proposition}{Proposition} [section]

\newtheorem{definition}{Definition} [section]

\newtheorem{remark}{Remark}[section]

\let\ssection=\section\renewcommand{\section}{\setcounter{equation}{0}\ssection}
\begin{document}
\date{}

\address{M. Darwich: Lebanese university, Beyrout Hadath.} \email{Mohamad.Darwich@lmpt.univ-tours.fr}
\title{Invariance of the White Noise for 
 the Ostrovsky  equation.}
\author{Darwich Mohamad.}
\begin{abstract}
In this paper, we construct
 invariant measures for the Ostrovsky equation associated with the norm $L^2$. On the other hand, 
we prove  the local well- posedness  in the besov space $\hat{b}^s_{p,\infty}$ for $sp >-1$.
 \end{abstract}

\maketitle

\section{Introduction}
 In this paper, we construct an invariant measure for a dynamical system defined by the Ostrovsky 
equation ($\text{Ost}$)

\begin{equation}\label{Ost}
\left\{\begin{array}{l}
\partial_tu - u_{xxx} + \partial_{x}^{-1}u +uu_x=0,\\
u(0,x)= u_0(x).
\end{array}\right.\end{equation}
for the quantities conserved by this equation .
The operator $ \partial_x^{-1}$ in the equation denotes a certain antiderivative with respect to the variable $x$ defined for $0$-mean value periodic
function the Fourier transform by 
$\widehat{(\partial_{x}^{-1}f)} = \frac{\hat f (\xi)}{i\xi}$.\\

Invariant measure play an important role in the theory of dynamical systems (DS). It is well known that the whole ergodic theory is based on this concept. On the other hand, they are necessary in various physical considerations.\\
In \cite{Zhidkov2} an infinite series of invariant measure associated with a higher conservation laws are
constructed for the one-dimensional Korteweg de Vries (KdV) equation:
$$
u_t+uu_x+u_xxx=0,
$$ by Zhidkov. In particular, invariant
measure associated to the conservation of the energie are constructed for this equation.\\
In the other hand, Tadahiro in \cite{OhKDV} was construct an invariant measure for (KdV) equation associated to the conservation of norm $L^2$ using the notion of the Wiener space. \\
 Equation \ref{Ost} is a perturbation of the  Korteweg de Vries (KdV) equation
with a nonlocal term and was deducted by Ostrovskii \cite{Ostrovskii} as a
model for weakly nonlinear long waves, in a rotating frame of reference, to describe the propagation of surface waves in the
ocean. \\

 The natural conserved quantities for (\ref{Ost}) are the :\\
$L^2$-norm:
$$ \|u(t)\|_{L^2}= (\int u^2(t,x) dx)^{\frac{1}{2}},$$
and the Hamiltonian:
$$H(u(t)) = \frac{1}{2} \int (u_x)^2 + \frac{1}{2}\int (\partial_x^{-1} u)^2 - \frac{1}{6}\int u^3.$$

We will construct invariant measures associated with the $L^2$-norm using the notion of Wiener Spaces.\\
The paper is organized as follows. In Section \ref{Notmesures} the
basic notation is introduced and the basic results are formulated. \\
 In section \ref{whitenoisewiener}
, we give the precise mathematical 
meaning of the white noise $Q$ defined for $ u = \sum_{n}a_ne^{inx}$ with mean $0$ by 
\begin{equation}\label{noise}
 dQ=\displaystyle{\frac{exp(-\frac{1}{2}\sum_{n\geqslant 1} \arrowvert a_n\arrowvert ^2)\Pi_{n\geqslant1}d a_n}
{\int exp(-\frac{1}{2}\sum_{n\geqslant1} \arrowvert a_n\arrowvert ^2)\Pi_{n\geqslant1}d a_n}},
\end{equation}
and show
that it is a (countably additive) probability measure on $\hat{b}_{p,\infty}^{s}$, $ p=2^{+}$, $s=-\frac{1}{2}+$ such that $sp<-1$, 
defined via the norm 
$$
\lVert f \lVert_{\hat{b}^s_{p,\infty}} := \lVert \hat{f} \lVert_{b^s_{p,\infty}} = \sup_{j} \lVert < n>^s \hat{f}(n)\|_{L^p_{|n| \sim 2^j}} = \sup_{j}
\bigg(\sum_{|n| \sim 2^j} <n>^{sp}|\hat{f}(n)|^p \bigg)^{\frac{1}{p}},
$$
and we go over the basic theory 
of abstract Wiener
spaces.
Finally in section \ref{localbsp}, we prove the  local well-posedness in 
$\hat{b}_{p,\infty}^{s}$.

\section{Notations and main results}\label{Notmesures}
We will use $C$ to denote various time independent constants, usually depending only upon $s$.
In case a constant depends upon other quantities, we will try to make it explicit.
 We use $A \lesssim B$ to denote an estimate of the form $A\leq C B$. similarly,
 we will write $A\sim B$ to mean $A \lesssim B$  and $B \lesssim A$. We  writre
 $\langle \cdot\rangle :=(1+|\cdot|^2)^{1/2}\sim 1+|\cdot| $. The notation $a^+$ denotes
$a+\epsilon$ for an arbitrarily small $\epsilon$. Similarly $a-$ denotes $a-\epsilon$.
On the circle, Fourier transform is
defined as
$$ \hat{f}(n)=\frac{1}{2\pi}\int_{\mathbb{T}}f(x)exp(-inx)dx$$
we  introduce the  Sobolev spaces   $H^{s}$  defined by~:
\begin{equation}\label{Esp1}
H^{s}=:\{ u \in \mathcal S^{'}(\mathbb{T}); ||u||_{H^{s}}
<+\infty \},
\end{equation}
where,
\begin{equation}\label{Ns}
||u|| _{H^{s}}= (2\pi)^{\frac{1}{2}}||{\langle.\rangle}^{s}\hat u||_{l^2_n}
\end{equation}
and 
$Y^{s} =:\{ u \in \mathcal S^{'}(\mathbb{T}); ||u||_{Y^{s}}
<+\infty \}$,
where
$$
||u|| _{Y^{s}}= ||u|| _{H^{s}} +||{\langle n\rangle}^{s}\hat u(n,\tau||_{l^2_nL^1_\tau},
$$

and the Besov-type space ${\hat{b}}^s_{p,\infty}$ defined via the norm
$$
\lVert f \lVert_{\hat{b}^s_{p,\infty}} := \lVert \hat{f} \lVert_{b^s_{p,\infty}} = \sup_{j} \lVert < n>^s \hat{f}(n)\|_{L^p_{|n| \sim 2^j}} = \sup_{j}
\bigg(\sum_{|n| \sim 2^j} <n>^{sp}|\hat{f}(n)|^p \bigg)^{\frac{1}{p}}.
$$
We will briefly remind the general construction of a Gaussian measure on a Hilbert space. Let $X$ be a Hilbert space, and $\{e_k\} $ be the orthonormal basis in $X$ which consists of eigenvectors of some operator $S=S^* > 0$ with corresponding  eigenvalues $0 < \lambda_1 \leq \lambda_2 \leq \lambda_3 ....\leq \lambda_k \leq ...$ We call a set $M \subset X$ a cylindrical
set iff:
$$
M =\{x \in X; [(x,e_1),(x,e_2),...(x,e_r)]\in F\}
$$
for some Borel $F \subset \mathbb{R}^{r}$, and some integer $r$. We define the measure $w$ as follows:
\begin{equation}\label{Gaussian}
w(M) = (2\pi)^{-\frac{r}{2}}\prod_{j=1}^{r}\lambda_j^{\frac{1}{2}}\int_{F}e^{-\frac{1}{2}\sum_{j=1}^{r}\lambda_{j}y_j^2}dy.
\end{equation}
One can easily verify that the class $\mathbb{A}$ of all cylindrical sets is an algebra on which the function $w$ is additive.
The function $w$ is called the centered Gaussian measure  on $X$ with the correlation operator $S^{-1}$.
\begin{definition}
The measure $w$ is called a  countably additive measure on an algebra $\mathbb{A}$ if $\lim_{n \rightarrow +\infty}(A_n)=0$ for any $A_n \in \mathbb{A}$($n=1,2,3...$) for which $A_1 \supset A_2 \supset A_3\supset.....\supset A_n\supset...$ and $\bigcap_{n=1}^{\infty}A_n = \phi$
\end{definition}
Now we give the following Lemma:
\begin{lemma}\label{lemmacountablyadditive}
The measure $w$ is countably additive on the algebra $\mathbb{A}$ iff $S^{-1}$ is an operator of trace class, i.e iff 
$\sum_{k=1}^{+\infty}\lambda_k^{-1} < +\infty $.
\end{lemma}
Now we present some definitions related to invariant  measure :
\begin{definition}
Let M be a complete separable metric space and let a function $h: \mathbb{R} \times M\longmapsto M$ for any fixed $t$ be a homeomorphism of the space $M$ into itself satisfying the properties:\\
\begin{enumerate}
\item $h(0,x)=x$ for any $x \in M$.
\item $h(t,h(\tau,x))=h(t+\tau,x)$ for any $t,\tau \in \mathbb{R}$ and $x \in M$.
\end{enumerate}
Then, we call the function $h$ a dynamical system with the space $M$. If $\mu$ is a Borel measure defined on the phase 
space $M$ and $\mu(\Omega)=\mu(h(\Omega,t))$for an arbitrary Borel set $\Omega \subset M$ and for all $t \in \mathbb{R}$, 
then it is called an invariant measure for the dynamical system $h$.
\end{definition}
\begin{definition}\label{Random}
  A random complex variable $Z$ is a measurable function from a probability space $(\Omega, F, P)$. to $\mathbb{C}$, i.e. it is a function 
$$
Z: \Omega \rightarrow \mathbb{C}
$$
such that for every Borel set $B \subset \mathbb{C}$, 
$$
Z^{-1}(B) =\{X \in B\}:=\{w \in \Omega: Z(w) \in B\} \in F.
$$
\end{definition}
\begin{definition}\label{distribution}
Let $Z$ be a random variable, we define the measure on Borel sets by 
$$
\mu_Z(B) = P\{Z \in B\} = P [Z^{-1}(B)].
$$
This measure is called the distribution of the random variable $Z$.
 \end{definition}
\begin{definition}\label{standard}
 Let $Z$ be a random variable, then $Z$ is called standard Gaussian radom variable on a probability space $(\Omega, F, P)$ if his probability desnsity function equals to 
$$
\frac{1}{2\pi}e^{-\frac{\mid z\mid ^2}{2}}.
$$
\end{definition}
Let us now state our results:

\begin{theorem}\label{thexistence1}
 Assume the mean $0$ condition on $u_0$. Let $ p=2^{+}$ and $ s > -\frac{1}{2}$ ($sp <-1$). Then, (\ref{Ost}) is locally well-posed in
 $\hat{b}_{p,\infty}^{s}$.
\end{theorem}
To prove This theorem, we can approximate the equation to the following problem:
\begin{equation}\label{Ostm}
\left\{\begin{array}{l}
\partial_tu^m - u^m_{xxx} + \partial_{x}^{-1}u^m + P_m(u^mu^m_x)=0,\\
u^m(0,x)=P_m u_0(x).
\end{array}\right.\end{equation}
Where $P_m$ is the projection onto the frequencies $|n| \leq N$ and $u^m = P_m u$. By Liouville's theorem, the finite dimensional white noise 
$$
d\mu_m = \displaystyle{C_m \exp (-\frac{1}{2} \int (u^m)^2 dx) }\displaystyle{\Pi_{x \in \mathbb{T}} du^m (x)}
$$
is invariant under of the flow of (\ref{Ostm}).
\begin{theorem}\label{theoremmeasure1}
Let $\{g_n(w)\}_n$ be a sequence of independent standard complex Gaussian random variables on a probability space $(\Omega, F, P)$.
Consider (\ref{Ost}) with initial data $u_0=\sum_{n\neq0}g_n(w)e^{inx}$ with mean $0$ condition, and $g_{-n}=\overline{g_n}$.
 Then, (\ref{Ost}) is globally well-posed almost surely in $w \in \Omega$. Moreover, the measure $Q$ defined in (\ref{noise}) is invariant under the flow.
 
\end{theorem}




\section{White noise measure and abstract wiener spaces}\label{whitenoisewiener}
In this section, we will construct an invariant measure related to the norm $L^2$. Let $Q$ denote the mean zero Gaussian white noise on the circle, i.e $Q$
 is the probability measure on distributions $u $with $\int_{\mathbb{T}} u = 0$ such that we have
$$
\int e^{i<f,u>}dQ(u)=e^{-\frac{1}{2}\rVert f \rVert_{L^2}^2}
$$
for all smooth mean $0$ function $f$ on $\mathbb{T}$. More precisely, the measure $Q$ is defined such that the fourier transform $\hat{Q}$ of $Q$ 
$$
\int e^{i<f,u>}dQ(u) = \hat{Q}(f) := e^{-\frac{1}{2}\rVert f \rVert_{L^2}^2}.
$$
Note that by the conservation of the mean, it follows that the solution $u(t)$ of (\ref{Ost}) has the spatial mean $0$ for all $t \in \mathbb{R}$
 as long as it exists.\\
 Let $u=\sum a_n e^{inx}$ be a real-valued function on $\mathbb{T}$ with mean $0$, i.e we have $a_0=0$ and 
$a_{-n}=\overline{a_n}$. Define $Q_N$ on $\mathbb{R}^{2N}$ with the density 
\begin{equation}\label{noise1}
  dQ_N =\displaystyle{\frac{exp(-\frac{1}{2}\sum_{n=1}^{N} \arrowvert a_n\arrowvert ^2)\Pi_{n=1}^Nd a_n}
{\int exp(-\frac{1}{2}\sum_{n=1}^{N} \arrowvert a_n\arrowvert ^2)\Pi_{n=1}^Nd a_n}}.
\end{equation}
Next, define the white noise $Q$ by 
\begin{equation}\label{noise2}
 dQ=\displaystyle{\frac{exp(-\frac{1}{2}\sum_{n\geqslant 1} \arrowvert a_n\arrowvert ^2)\Pi_{n\geqslant1}d a_n}
{\int exp(-\frac{1}{2}\sum_{n\geqslant1} \arrowvert a_n\arrowvert ^2)\Pi_{n\geqslant1}d a_n}}.
\end{equation}

Let $\dot{H}_0^s$ be the homogeneous Sobolev space restricted to the real-valued mean $0$ elements. Let $<.,.>_{\dot{H}_0^s}$ denote the inner 
product in $\dot{H}_0^s$ i.e $<\sum c_ne^{inx},\sum d_n e^{inx}>_{\dot{H}_0^s}= \sum \mid n \mid^{2s} c_n\overline{d_n}$.\\
Let $B_s = \sqrt {-\Delta}^{2s}$, then the weighted exponentials $\{\mid n \mid^{-s} e^{inx}\}_{n\neq0}$ are the eigenvectors of $B_s$ with the eigenvalue 
$\mid n \mid ^{2s}$, forming an orthonormal basis of $\dot{H}_0^s$. Note that
$$
-\frac{1}{2}<B^{-1}\phi,\phi>_{\dot{H}_0^s} = -\frac{1}{2}\sum \mid a_n\mid^2.
$$
The right hand side is exactly the expression appearing in the exponent in (\ref{noise2}). By Lemma \ref{lemmacountablyadditive} 
$Q$ is countably additive if and only if $B$ is of trace class, i.e $\sum \mid n \mid^{2s} <\infty$.
 
Hence, $\text{supp}(Q)\subset\bigcap_{s<-\frac{1}{2}}H^s$ is a natural space to work on, and
it is known (cf Zhidkov \cite{Zhidkov}) that the white noise $Q$ is supported on $ \bigcap_{s<-\frac{1}{2}}H^s$.
Unfortunately, we cannot prove 
a local-in -time solution of (\ref{Ost}) in $H^{s}$, $s<-\frac{1}{2}$ (see \cite{Darwichostr}).\\
Then for this and Remark \ref{distributionofu} also, following Oh \cite{OhKDV}, we propose to work on a Besov-type space ${\hat{b}}^s_{p,\infty}$ defined via the norm
$$
\lVert f \lVert_{\hat{b}^s_{p,\infty}} := \lVert \hat{f} \lVert_{b^s_{p,\infty}} = \sup_{j} \lVert < n>^s \hat{f}(n)\|_{L^p_{|n| \sim 2^j}} = \sup_{j}
\bigg(\sum_{|n| \sim 2^j} <n>^{sp}|\hat{f}(n)|^p \bigg)^{\frac{1}{p}}.
$$

By Hausdorff-Young's inequality, we have $ \hat{b}^s_{p,\infty} \supset B^s_{p^{\prime},\infty} $ for $p > 2$, where $B^s_{p^{\prime},\infty}$ 
is the usual Besov space with $p^{\prime}= \frac{p}{p-1}$ .\\

 This space has many advantages, it contains the support of the white noise for $sp < -1$.This follows from the 
theory of abstract Wiener spaces (c.f Gross, Kuo).\\
Since ${\hat{b}}^s_{p,\infty}$ is not a Hilbert space, we need to go over the basic theory of abstract Wiener
spaces.
\begin{remark}\label{lemmenonexistence}
In view of the results in \cite{Darwichostr} , we can not hope to have a local-time solution via the fixed point argument in $\mathbb{H}^{s}$, $ s<-\frac{1}{2}$. 
\end{remark}
\begin{remark}\label{distributionofu}
 Note that the measure $Q$ can be defined as  probability distribution  for the following random variable
$$
u: \Omega \times \mathbb{T} \longmapsto \mathbb{R}\\
      (w,x)                 \longmapsto \sum_{n} g_n(w) e^{inx},
$$

where $g_n$ are independent standard complex Gaussian random variables.
By another way, remark that 
$E(\Arrowvert u \Arrowvert_{H^s}^2)= E(\sum _{n}\frac{\arrowvert g_n\arrowvert^2}{\arrowvert n \arrowvert^{-2s}}) \lesssim 
\sum _{n} \frac{1}{\arrowvert n \arrowvert^{-2s}} < +\infty$, iff $s < -\frac{1}{2}$, but by Remark \ref{lemmenonexistence} we cannot work in this space.	
On the other side
$E(\Arrowvert u \Arrowvert_{b^s_{p,\infty}}^p)= E(\sum _{n}\frac{\arrowvert g_n\arrowvert^p}{\arrowvert n \arrowvert^{-ps}}) \lesssim 
\sum _{n} \frac{1}{\arrowvert n \arrowvert^{-ps}} < +\infty$, iff $sp < -1$. Indeed for $ sp<-1$ the map $\sum_n g_n(w) e^{inx}$
defines a ( Gaussian ) measure on $b^s_{p,\infty}$ , i.e the white measure $Q$.
 
\end{remark}
 
\subsection{Global existence in $\hat{b}_{p,\infty}^{s}$}
Now we will prove the global existence using the local existence (proved in section 6) and the invariance of the white noise.
\begin{remark}
 Note that by Liouville's theorem, the Lebesgue measure $\Pi_{n\geqslant1}^{N}d a_n$ is invariant under the flow of (\ref{Ostm}). Hence,
 the finite dimensional version $Q_N$ of $Q$ is invariant under the flow of (\ref{Ostm}). 
\end{remark}
Using the invariance of $Q_N$, we have the 
following estimate on $u^N$
\begin{proposition}\label{Growth}
 Given $T > 0$ and $\epsilon > 0$, there exists $\Omega_N \subset \hat{b}_{p,\infty}^{s} $ with $Q_n(\Omega_N^c)<\epsilon$ such that for $u_0^N \in \Omega_N$, (\ref{Ostm}) is well-posed on
$[-T,T]$, with the following growth:
$$
\parallel u^N(t)\parallel_{\hat{b}_{p,\infty}^{s}} \lesssim (log\frac{T}{\epsilon})^{\frac{1}{2}},~~ for \mid t \mid \leq T.
$$
\end{proposition}
In proving Proposition \ref{Growth}, we need to assume the following estimate:
\begin{lemma}\label{Growth1}
There exists $c > 0$, independent of $N$, such that for sufficiently large $K > 0$
$$
Q_N(\{\parallel u_0^N\parallel_{\hat{b}_{p,\infty}^{s}}  > K\}) < e^{-cK^2}.
$$
\end{lemma}
\proof See \cite{Fernique}.\\
\textbf{Proof of proposition \ref{Growth}}. Let $S_N(t)$ denote the flow map of (\ref{Ostm}), and define \\
$$
\displaystyle{\Omega_N = \bigcap _{j=-[T\setminus\delta]}^{[T\setminus\delta]} S_N^j(\delta)(\{\parallel u_0^N\parallel_{\hat{b}_{p,\infty}^{s}} \leq K\})}.
$$
By invariance of $Q_N$ and $\delta \sim K^{-\alpha}$, we have 
$$
Q_N(\Omega_N ^c) \lesssim \frac{T}{\delta}Q_N(\{\parallel u_0^N\parallel_{\hat{b}_{p,\infty}^{s}}  > K\}) \sim T K^{\alpha}e^{-cK^2}.
$$
By choosing $K \sim (\text{log}\frac{T}{\epsilon})^{\frac{1}{2}}$, we have $Q_N(\Omega_N ^c) < \epsilon$. Moreover,
 by its construction, $\rVert u^N(j\delta)\rVert_{\hat{b}_{p,\infty}^{s}} \leq K$ for $j= 0,..., \pm [T\setminus \delta]$. By local theory, we have
$$
\parallel u^N(t)\parallel_{\hat{b}_{p,\infty}^{s}} \leq 2K \sim (\text{log}\frac{T}{\epsilon})^{\frac{1}{2}} ~~ for |t| \leq T,
$$ 
then $ \Omega_N$ has the desired property.\\
As a corollary to Proposition \ref{Growth}, one needs to prove the following statements.
a) For $\epsilon >0$, there exists $\Omega_\epsilon \subset \hat{b}_{p,\infty}^{s} $ with $Q(\Omega_\epsilon^c) < \epsilon$ such that for any $u_0 \in 
\Omega_\epsilon$, (\ref{Ost}) is globaly well-posed with the growth estimate:
\begin{equation}\label{Growth2}
\parallel u(t) \parallel_{\hat{b}_{p,\infty}^{s}} \lesssim \bigg( \text{log} \frac{1 + \mid t \mid}{\epsilon}\bigg)^{\frac{1}{2}}, ~~
\text{for all} ~ t \in \mathbb{R} 
 \end{equation}
b) The uniform convergence lemma:
$$
\parallel u - u^N\parallel_{C([-T,T];\hat{b}_{p,\infty}^{s})} \longrightarrow 0
$$
as $N \longrightarrow +\infty$ uniformly for $u_0 \in \Omega_\epsilon$.\\

Note that (a) implies that the problem is globaly well posed, since $\cup_{\epsilon >0} \Omega_\epsilon$ has probability 1. We can prove (a) and (b) by 
estimating the difference $u-u^N$, using the linear and bilinear estimates and applying Proposition \ref{Growth} to $u^N$. Note that the nonlinearity of the 
difference equation is given by 
$$
R(t)=\partial_x u^2(t) -P_N(\partial_x(u_N)^2(t)).
$$
Since $P_N\bigg((P_N u)^2\bigg) = (P_\frac{N}{2}u)^2$, we have
$$
R(t)= \partial_x(u^2-(P_\frac{N}{2}u)^2) + P_N\partial_x\bigg((P_\frac{N}{2}u)^2 - u^2\bigg) + P_N\partial_x(u^2 - (u^N)^2).
$$
After applying the nonlinear estimate, the first two terms can be made small due to the factor $u-P_{\frac{N}{2}}u$, 
and the last term has the factor $u-u^N$,
which we need to close the argument.\\


\subsection{Abstract Wiener spaces}\label{Wiener}
In section 2, we reviewed the Gaussian measures in
Hilbert spaces. However, ${\hat{b}}^s_{p,\infty}$ is not a Hilbert space, so we briefly go over the basic theory
of abstract Wiener spaces.
Recall the following definitions from Kuo: Given a real separable Hilbert space $H$ with norm 
$\rVert.\rVert $, let $\mathbb{F}$ denote the set of finite dimensional orthogonal projections $\mathbb{P}$ of H.
Then, define a cylinder set $E$ by $E = \{x \in H : \mathbb{P}x \in \mathbb{F} \}$ where $\mathbb{P} \in \mathbb{F}$ and $\mathbb{F}$ is a Borel
subset of $\mathbb{P}H$, and let $R$ denote the collection of such cylinder sets. Note that $R$ is a field
but not a $\sigma$-field. Then, the Gauss measure $\mu$ on $H$ is defined by
$$
\mu(E)= (2\pi)^{-\frac{n}{2}}\int_F e^{-\frac{\lVert x\lVert}{2}}dx
$$
for $E \in R$, where $n = \text{dim} \mathbb{P}H$ and $dx$ is the Lebesgue measure on $\mathbb{P}H$. It is known that $μ$
is finitely additive but not countably additive in $R$.
A seminorm $|||.|||$ in $H$ is called measurable if for every $\epsilon > 0$, there exists $\mathbb{P}_\epsilon \in F $ such
that$$
\mu(|||Px||| > \epsilon) < \epsilon$$
for $\mathbb{P} \in F$ orthogonal to $\mathbb{P}_{\epsilon}$ . Any measurable seminorm is weaker than the norm of $H$,
and $H$ is not complete with respect to$ ||| .|||$ unless $H$ is finite dimensional. Let $B$ be the
completion of $H$ with respect to $|||.|||$ and denote by $i$ the inclusion map of $H$ into $B$. The
triple $(i, H, B)$ is called an abstract Wiener space.\\
Now, regarding $y \in B ^{∗}$ as an element of $H ^{*} \equiv H$ by restriction, we embed $B^{*}$ in $H$.
Define the extension of $\mu$ onto $B$ (which we still denote by $\mu$) as follows. For a Borel set
$F \subset \mathbb{R}^n$ , set
$$
\mu(\{x \in B : ((x, y_1 ), ... , (x, y_n )) \in F \}) := \mu(\{x \in H : <x, y_1>_H, · · ·
, <x, y_n>_H)
\in F \})$$
where $y_j$ 's are in $B^*$ and $(.,.)$ denote the natural pairing between $B$ and $B^*$.  Let $R_B$ denote
the collection of cylinder sets $\{x \in B : ((x, y_1 ), ... , (x, y_n )) \in F \}$ in $B$. Note that the pair
$(B, \mu)$ is often referred to as an abstract Wiener space as well.
\begin{theorem}\label{countablyadditive}(Gross \cite{Gross}).
The measure $Q$ defined in (\ref{noise2}) is countably additive in the $\sigma$-field generated by $R_B$.
 \end{theorem}
From now, let $H= L^2(\mathbb{T})$ and $B={\hat{b}}^s_{p,\infty}(\mathbb{T})$ for $sp < -1$.
\begin{proposition}\label{normmeasurable}
 The semi norm $\|.\|_{\hat{b}^s_{p,\infty}}$ is measurable for $sp< -1$.
\end{proposition}
Hence, $(i, H, B) = (i,L^2,\hat{b}^s_{p,\infty})$ is an abstract Wiener space, and $Q$ defined in (\ref{noise2}) is
countably additive in $\hat{b}^s_{p,\infty}$ . We present the proof of Proposition \ref{normmeasurable} at the end of this
subsection. For our application, we can choose $s$ and $p$ such that $sp < -1$. Note that follows from the proof that $(i,L^2,\hat{b}^s_{p,\infty})$ 
is an abstract Wiener space for $sp < -1$.
\\

\begin{theorem}
 Let $(i, H, B)$ be a Wiener space. Then,
there exists $c > 0$ such that $\int_B e^{c\| x\|_B^2} Q(dx) < \infty$. Hence, there exists $A > 0 $such that
$Q(\|x\|_B> K) \leq e^{−A K^2}$ for sufficiently large  $A> 0$.

\end{theorem}
\proof See Theorem 3.1 in \cite{Kuo}.\\
 To prove Proposition \ref{normmeasurable}, we  need the followings lemma:
\begin{lemma}\label{lemmagn}(Lemma 4.7 in \cite{Ohshrodinger})
 Let $\phi = \sum_n g_n e^{inx}$, where $(g_n)_{n=1}^{\infty}$ is a sequence of independent standard complex-valued Gaussian random variables. Then, for $M$ 
dyadic and $\delta >0$, we have 
$$
\lim _{M \longrightarrow \infty} M^{1-\delta} \frac{max_{|n| \sim M}|g_n|^2}{\sum_{|n| \sim M} |g_n|^2} = 0,~~ a.s.
$$
\end{lemma}
Now, we present a large deviation lemma,. This can be proved by a direct computation using the polar coordinate . See \cite{Bourgain}.
\begin{lemma}\label{lemmagn1}
 Let $M$ be a dyadic, and $R =R(M) \geq M^{\frac{1}{2}+}$. Then, there exists $c$ such that 
$$
\mathbb{P}_w\bigg[\bigg(\sum_{n\sim M} |g_n(w)|^2\bigg)^{\frac{1}{2}} \geq R \bigg]\leq e^{-cR^2}
$$
for all dyadic $M$ (i.e $c$ independent of $M$). Moreover, this is essentially sharp in the sense that the estimation can not hold if
 $R \leq M^{\frac{1}{2}}$.
\end{lemma}
\textbf{Proof of Proposition \ref{normmeasurable}}: Let $\epsilon >0$,  it suffices to show that for given $\epsilon >0$, there exists large $M_0$ such that
$$
Q\bigg (\| \mathbb{P}_{>M_0}\phi\|_{\hat{b}^s_{p,\infty}} > \epsilon\bigg) <\epsilon,
$$
where ${P}_{>M_0}$ is the projection onto the frequencies $|n| > M_0$. By Egoroff's theorem ( cited in Section 2), there exists a set $E$ such that 
$Q(E^c) < \frac{1}{2}\epsilon$.
Fix $K > 1$ and $\delta \in (0,\frac{1}{2})$, then by Lemma \ref{lemmagn} there exists $M_0$ large enough such that
\begin{equation}\label{proofprop0}
\frac{\|\{g_n(w)\}_{|n| \sim M}\|_{L^\infty_n}}{\|\{g_n(w)\}_{|n| \sim M}\|_{L^2_n}} \leq M^{1-\delta},
\end{equation}
for all $w \in E$ and dyadic $M> M_0$.\\
The basic idea of the following argument is due to Bourgain's dyadic pigeonhole principle in \cite{Bourgain}.  Let $\{\sigma_j\}_{j\geq1}$ 
be a sequence of positive numbers such that $ \sum \sigma_j = 1$, and let $M_j = M_02^j$. Note that $\sigma_j= C2^{-\lambda j}= CM_0^{\lambda}M_j^{-\lambda}$
for some small $\lambda >0$ ( to be determined later). Then, we have
\begin{align}\label{proofprop}
 Q\bigg( \| \mathbb{P}_{>M_0}\phi\|_{\hat{b}^s_{p,\infty}} > \epsilon\bigg) &\leq Q\bigg(\|\{g_n\}_{|n| >M_0}\|_{\hat{b}^s_{p,1}} > \epsilon\bigg)\nonumber\\
&\leq \sum_{j=0}^{\infty} Q(\|\{<n>^sg_n\}_{|n|\sim M_j}\|_{L^p_n} > \sigma_j\epsilon\bigg),
\end{align}
where ${\hat{b}^s_{p,1}}$ is as ${\hat{b}^s_{p,\infty}}$ with the $l^\infty$ norm over the dyadic blocks replaced by the $l^1$ sum.\\
By interpolation and (\ref{proofprop0}), we have 
\begin{align}
 \|\{<n>^sg_n\}_{|n|\sim M_j}\|_{L^p_n} &\sim M_j^s \| \{g_n\}_{|n|\sim M_j}\|_{L^p_n} \leq M_j^s\| \{g_n\}_{|n|\sim M_j}\|^{\frac{2}{p}}_{L^2_n}
\| \{g_n\}_{|n|\sim M_j}\|_{L^\infty_n}^{\frac{p-2}{p}}\nonumber\\
&\leq M_j^s \| \{g_n\}_{|n|\sim M_j}\|_{L^2_n}
\bigg( \frac{\|{g_n(w)}_{|n| \sim M}\|_{L^\infty_n}}{\|\{g_n(w)\}_{|n| \sim M}\|_{L^2_n}}\bigg)^{\frac{p-2}{p}} \leq M_j^{s-\delta\frac{p-2}{p}}
\| \{g_n\}_{|n|\sim M_j}\|_{L^2_n}\nonumber\\
\end{align}
a.s thus, if $\|\{<n>^sg_n\}_{|n|\sim M_j}\|_{L^p_n} > \sigma_j \epsilon$, we obtain that $\| \{g_n\}_{|n|\sim M_j}\|_{L^2_n} \gtrsim R_j$ where
$R_j:= \sigma_j \epsilon M_j^{-s + \delta\frac{p-2}{p}}$. For $p= 2 + 2\theta$, we have  $-s + \delta\frac{p-2}{p} = \frac{-sp + 2\delta\theta}{2+2\theta} > 
\frac{1}{2}$ by taking $\delta$ sufficiently close to $\frac{1}{2}$ since $-sp > 1$. By taking $ \lambda >0$ sufficiently small, 
$R_j = \sigma_j \epsilon M_j^{-s + \delta\frac{p-2}{p}} = C\epsilon M_0^\lambda M_j^{-s + \delta\frac{p-2}{p}-\lambda}\gtrsim 
C\epsilon M_0^\lambda M_j^{\frac{1}{2}+}$. Then, by Lemma \ref{lemmagn1}, we have
\begin{equation}\label{proofprop1}
 Q\bigg( \| \mathbb{P}_{>M_0}\phi\|_{\hat{b}^s_{p,\infty}} > \epsilon\bigg) \leq \sum_j^{\infty} e^{-C^2M_0^{1+2\lambda}+2^j+\epsilon^2} \leq \frac{1}{2}\epsilon,
\end{equation}
by choosing $M_0$ sufficiently large.\\

 \textbf{Proof of Theorem \ref{theoremmeasure1}}:

By the invariance of $Q_N$ and the uniform convergence of $u^N$ to $u$, we obtain the invariance of $Q$

\section{Local well-Posedness in $\hat{b}_{p,\infty}^{s}$}\label{localbsp}

In this section, we prove Theorem \ref{thexistence1} via the fixed point argument.  We
go over the previous local well-posedness theory of Ostrovsky equation to motivate the definition of the
Bourgain space $W_p^{s,b}$ with the weight, adjusted to $\hat{b}_{p,\infty}^{s}$ .\\
 
{
We have proved that for $n=n_1+n_2$ : $$\sigma(\tau,\tau_1,n,n_1) =
\text{max}\{\mid \tau + m(n)\mid, \mid \tau_1 + m(n_1)\mid, \mid \tau - \tau_1 +  m(n-n_1)\mid\} 
\gtrsim \mid nn_1n_2\mid.$$
Recall that this estimates implies that:
\begin{equation}\label{optimal}
\frac{\mid n \mid^{s+1}\mid n_1(n-n_1)\mid^{-s}}{< \tau + m(n) >^{\frac{1}{2}}< \tau_1 + m(n_1) >^{\frac{1}{2}} 
<\tau- \tau_1 + m(n_2) >^{\frac{1}{2}}} \lesssim \frac{\mid n \mid< n >^s}{< n_1 >^s< n_2 >^s}\frac{1}{\sigma^{\frac{1}{2}}} \lesssim 1
\end{equation}
for $s \geq -\frac{1}{2}$.  Note that (\ref{optimal}) is optimal, for example, when $<\tau + m(n) >\sim <3nn_1n_2>$ and 
 $< \tau_j + m(n_j) > \ll <3nn_1n_2>^{0^+}$. To exploit this along with the fact the free solution concentrates on the curve $\{ \tau =n^3 - \frac{1}{n}\}$,
we define the weight $v(n,\tau)$ by :
$$
v(n,\tau) = 1+ \sum_{k \neq 0} min(<k>,<n-k>)^{\delta}1_{A_k},
$$
where 
$$
A_k=\{(n,\tau):\mid n \mid \geq 1, <\tau-n^3 +\frac{1}{n} +3n(n-k)k> \ll <n>^{\frac{1}{100}}\},
$$
and $\delta = 0^{+}$ (to be determined later).

Note that, for fixed $n$ and $\tau$, there are at most two values of $k$ such that $|(n-k)k + \frac{\tau-n^3 + \frac{1}{n}}{3n} | \ll <n>^{-1 +\frac{1}{100}}$.\\
It follows from the definition  that $v(n,\tau) \lesssim \text{max}(1,(\frac{<\tau-n^3+\frac{1}{n}}{<n>})^{0^+})\leq < \tau - n^3 +\frac{1}{n}>^{0^+}$.\\
Now, define the Bourgain space $W_p^{s,b}$ with the weignt $v$ via the norm:
$$
\parallel u \parallel_{W_{p}^{s,b}} = \parallel \hat{u}\parallel_{\hat{W}_p^{s,b}} := \parallel v\hat{u}\parallel_{\hat{X}_p^{s,b}} 
+ \parallel\hat{u}\parallel_{\hat{Y}_p^{s,b-\frac{1}{2}}},
$$
where 
$$\parallel f\parallel_{\hat{X}_p^{s,b}} := \sup_{j}\parallel<n>^s<\tau-n^3 + \frac{1}{n}>^bf(n,\tau)\parallel_{L^p_{\mid n \mid \sim 2^j}L^p_\tau}
$$
and 
$$\parallel f\parallel_{\hat{Y}_p^{s,b}} := \sup_{j}\parallel<n>^s<\tau-n^3 + \frac{1}{n}>^bf(n,\tau)\parallel_{L^p_{\mid n \mid \sim 2^j}L^1_\tau}.
$$
For our paper, we set $b=\frac{1}{2}$. Note that $\hat{Y}_p^{s,0}$ is introduced so that we have $W^{s,\frac{1}{2}}_{p}(\mathbb{T}\times [-T,T]) \subset 
C([-T,T];\hat{b}_{p,\infty}^{s}(\mathbb{T}))$. 
\subsection{Linear Estimates} Let $S(t)$ the free evolution of Equation \ref{Ost} and $\eta(t)$ be a smooth cutoff such that $\eta(t)=1$ on 
$[-\frac{1}{3},+\frac{1}{3}]$ and =0 for $|t| \geq 1$.\\
\begin{lemma}
 For any $s \in \mathbb{R}$, we have $\parallel \eta S(t)u_0\parallel_{{W_{p}^{s,\frac{1}{2}}}}\lesssim \parallel u_0\parallel_{\hat{b}_{p,\infty}^{s}}$.
\end{lemma}
\proof: Recall that $ v(n,\tau) \lesssim <\tau - n^3 + \frac{1}{n}>^{0^+}$. Noting that $\hat{(\eta(t) S(t)u_0)}(n,\tau) = \hat{\eta}(\tau-m(n))\hat{u_0}(n)$,
we have 
\begin{align}
 \parallel \eta(t)S(t)u_0\parallel_{{W_{p}^{s,\frac{1}{2}}}} &\leq \sup_{j}\parallel <n>^s<\tau-n^3 + \frac{1}{n}>^{\frac{1}{2}^+}
\hat{\eta}(\tau-m(n))\parallel_{L^p_\tau}\mid\hat{u_0}(n)\mid\parallel_{L^p_{\mid n\mid \sim 2^j}}\nonumber\\
&+\sup_j \parallel<n>^s\parallel\hat{\eta}(\tau-n^3 + \frac{1}{n})\parallel_{L^1_\tau}\mid\hat{u_0}(n)\mid\parallel_{{L^p}_{\mid n\mid\sim 2^j}} 
\leq C_\eta \parallel u_0\parallel_{\hat{b}_{p,\infty}^{s}},\nonumber
\end{align}
where $ C_\eta = \parallel < \tau>^{\frac{1}{2}^+}\hat{\eta}(\tau)\parallel_{L^p_\tau} + \parallel \hat{\eta}\parallel_{L^1} < \infty$.
\\
Now, we estimate the Duhamel term. By the standard computation, we have
\begin{align}
\int_{0}^{t}S(t-t^{\prime})dt^{\prime} &=-i \sum_{k=1}^{\infty}\frac{i^kt^k}{k!}\sum_{n\neq0}e^{i(nx + m(n)t})\int \eta(\lambda - m(n))\hat{F}(n,\lambda)d\lambda\nonumber\\
&+i\sum_{n\neq0}e^{inx}\int\frac{(1-\eta)(\lambda-m(n))}{\lambda - m(n)}e^{i\lambda t}\hat{F}(n,\lambda)d\lambda\nonumber\\
&+i\sum_{n\neq0}e^{i(nx + m(n)t)}\int\frac{(1-\eta)(\lambda-m(n))}{\lambda - m(n)}e^{i\lambda t}\hat{F}(n,\lambda)d\lambda\nonumber\\
&=: F_1 + F_2 + F_3.\nonumber
\end{align}
\begin{lemma}
 For any $ s \in \mathbb{R}$, we have
$$
\parallel \eta(t)F_1\parallel_{{W}_p^{s,\frac{1}{2}}}, \parallel \eta(t)F_2\parallel_{{W}_p^{s,\frac{1}{2}}},
 \parallel \eta(t)F_3\parallel_{{W}_p^{s,\frac{1}{2}}} \lesssim \parallel \eta(t)F\parallel_{{W}_p^{s,-\frac{1}{2}}}
$$
\end{lemma}
\proof See \cite{Colliander}.
\subsection{Bilinear estimate}
In this part, we will treat the bilinear term $\partial_x(uv)$ and we have the following one:
\begin{proposition}Let $u_1$ and $u_2$ have the spatial means $0$ for all $ t$ in $\mathbb{R}$. Then, there exist $s = -\frac{1}{2} +$, $p > 2$ with $sp<-1$, and
$\theta > 0$ such that
\begin{equation}\label{Bilinearlp}
\parallel \eta _{2T} \partial_x(u_1u_2)\parallel_{{W}_p^{s,-\frac{1}{2}}} 
\lesssim T^{\theta}\parallel u_1\parallel_{{W}_p^{s,\frac{1}{2}}}\parallel u_2\parallel_{{W}_p^{s,\frac{1}{2}}}. 
 \end{equation}
\end{proposition}
To prove this Proposition, we need the followings lemmas:
\begin{lemma}\label{int}(Ginibre-Tsutsumi-Velo \cite{Ginibre})
 Let $ 0 \leq \alpha \leq \beta$ and $\alpha + \beta > \frac{1}{2}$. Then, we have
$$
\int <\tau>^{-2\alpha}<\tau-a>^{-2\beta}d\tau \lesssim~~ <a>^{-\gamma},
$$
where $ \gamma = 2\alpha-[1-2\beta]_{+}$ with $[x]_+=x$ if $x>0$, $= \epsilon >0$ if $x=0$, and $=0$ if $x<0$.
\end{lemma}
\begin{lemma}\label{sum}
 For $l_1 + 2l_2 >1$ with $l_1,l_2 >0$, there exists $c>0$ such that for all $n\neq 0$ and $\lambda \in \mathbb{R}$, we have
\begin{equation}\label{lemmal1l2}
\sum_{n_1\neq0,n}\frac{1}{<n_1>^{l_1}<\lambda+n_1(n-n_1)>^{l_2}} \leq c.
\end{equation}
\end{lemma}
\proof: If $l_1 > 1$, then (\ref{lemmal1l2}) is clear. If $l_2 > \frac{1}{2}$, let $\alpha_1, \alpha_2$ the roots of the polynomial
$$
\lambda + n_1(n-n_1) = 0.
$$
There are at most 10 $n_1's$ such that $\mid n_1 - \alpha \mid \leq 2$ or $ \mid n_1 - \beta \mid \leq 2$. The remaining $n_1's$  satisfy
$$
(1+ \mid (n_1 - \alpha)(n_1 - \beta)\mid) \geqslant \frac{1}{2}(1 + \mid n_1 + \alpha \mid)(1 + \mid n_1 + \beta \mid ).
$$
Now we have\\ 
\begin{align}
(\ref{lemmal1l2})  \lesssim \sum_{n_1} 
\frac{1}{<\lambda+ n_1(n-n_1)>^{l_2}} &= \sum_{n_1} \frac{1}{\big(1+ \mid (n_1 - \alpha)(n_1 - \beta)\mid)\big)^{l_2}}\nonumber\\
&\lesssim 
\sum_{n_1}\frac{1}{(1 + \mid n_1 + \alpha \mid)(1 + \mid n_1 + \beta \mid )^{l_2}}.\nonumber
\end{align}
Hence, appying the Cauchy-Schwarz inequality we obtain the desired result.\\
Now assume that $ l_1 \in (0,1]$ and $l_2 \in (0,\frac{1}{2}]$, since $l_1 +2l_2 > 1$, there exists $\epsilon > 0$ such that $l_1 +2l_2 -3\epsilon \geqslant 1$.\\
If $ P_{n,\tau}(n_1):= \tau +n_1(n-n_1)$ has two real roots, i.e $P_{n,\tau}(n_1) = -(n_1-r_1)(n_1-r_2)$, then there are at most 6 values of $n_1$ such that
$\mid n_1 - r_j\mid \leq 1$. For the remaining values of $n_1$, we have $<P_{n,\tau}(n_1)> > \frac{1}{4} \prod_{j=1}^2 <n_1 - r_j>$. Then, (\ref{lemmal1l2})
follows from Holder inequality with $p=(l_1-\epsilon)^{-1}$ and $ q= (l_2 -\epsilon)^{-1}$, we have 
$$
\text{LHS of }~ (\ref{lemmal1l2}) \lesssim \bigg(\sum_{n_1} <n_1>^{-pl_1}\bigg)^{\frac{1}{p}}\prod_{j=1}^{2}
\bigg(\sum_{n_1} < (n_1-r_j>^{-ql_2}\bigg)^{\frac{1}{q}} < \infty
$$
since $ pl_1 > 1$ and $ ql_2 > 1$.\\
If $P_{n,\tau}(n_1)$ has only one or no real root, then we have $|P_{n,\tau}(n_1)| \geqslant (n_1 - \frac{1}{2} n)^2$ for all $n_1 \in \mathbb{Z}$. Then, by 
Holder inequality with $p=(l_1-\epsilon)^{-1}$ and $ q=(2l_2 - 2\epsilon )^{-1}$, we have
$$
\text{LHS of} ~(\ref{lemmal1l2}) \leq \bigg(\sum_{n_1} <n_1>^{-pl_1}\bigg)^{\frac{1}{p}}(\sum_{n_1} < (n_1-\frac{1}{2}n)^2>^{-ql_2}\bigg)^{\frac{1}{q}} < \infty
$$
since $pl_1 > 1$ and $2ql_2 = \frac{l_2}{l_2 - \epsilon} > 1$.

\begin{lemma}\label{intOmega}
 Let $\textit{O(n)}= \{\varsigma \in \mathbb{R}: \varsigma=-3nn_1n_2 + o(<nn_1n_2>^{\frac{1}{100}}) , \text{for some} ~ n_1 \in \mathbb{Z}\text{ with}~ n = n_1+n_2 \}$, then
$$\int <\tau -n^3 +\frac{1}{n}> ^{-\zeta} 1_{\textit{O(n)}}(\tau -n^3 +\frac{1}{n})d\tau \lesssim 1, \zeta = 1-.
$$
\end{lemma}
\proof: Let
$$ \textit{K(n)}= \{\varsigma \in \mathbb{R}, \left|\varsigma\right|\sim M, \varsigma =-3n_1n_2n + O(<nn_1n_2>^{\frac{1}{100}}\},$$
where $M \geq 1 $ is a dyadic number.\\
We will prove that the measure of $\textit{K(n)} \leq M^{1-o}$, for some $o$ small. The hypotheses are symmetric in $n_1$, $n_2$ so we may assume $\left| n_1\right| \geq \left| n_2\right|$. 
First suppose that $|n_1| \leq |n|$.

Since $n_1$, $n_2
\in \mathbb{Z}$ and $\left|n_1n_2n\right| \lesssim \left|n\right|^3$ we conclude that $\left|n\right| \lesssim \left|\varsigma\right|\lesssim \left|n\right|^3$.\\
Suppose $\left|\varsigma\right|\sim M$ and $\left|n\right|\sim N$(dyadic). We have that $M\sim N^p$ for some $p \in [1,3]$, then $\left|n_1n_2\right|\sim M^{1-\frac{1}{p}}$.
We make the crude observation that there are at most $M^{1-\frac{1}{p}}$ multiplies $M^{\frac{1}{p}}$ in the dyadic block $\{\left|\varsigma\right|\sim M\}$.
Hence, the set of possible $\varsigma =-3nn_1n_2 + o(<nn_1n_2>^{\frac{1}{100})}$ must lie inside a union of $M^{1-\frac{1}{p}}$intervals of size $M^{\frac{1}{100}}$, each
of which contains an integer multiple of $n$. We have then that 
$$
\left| \textit{K(n)}\right| < M^{1-\frac{1}{p}}M^{\frac{1}{100}} \lesssim M^{\frac{3}{4}}.
$$
In case $\left|n\right| \leq \left|n_1\right|$, we must have $\left|n_1\right|\lesssim \left|\varsigma\right| \lesssim \left|n_1\right|^3$ so, if $\left| n_1\right| \sim N_1$(dyadic), we must have $M \sim N_1^p$ for some $p \in [1,3]$ and we can repeat the argument presented above.\\
Now by  a changement of variables \\
$\int <\tau -n^3 +\frac{1}{n}> ^{-\zeta} 1_{\textit{O(n)}}(\tau -n^3 +\frac{1}{n})d\tau = \int <z>^{-\zeta} 1_{\textit{O(n)}}(z)dz = \displaystyle{\sum
_{M dyadic}} \int_{\left|z\right|\sim M}<z>^{-\zeta} 1_{\textit{O(n)}}(z)dz \leq \displaystyle{\sum_{M dyadic}} M^{-1}M^{1-o} \lesssim 1$.\\
\textbf{Proof of Proposition \ref{Bilinearlp}:} First, we prove
\begin{equation}\label{Bilinearlp1}
\parallel  \partial_x(u_1u_2)\parallel_{{W}_p^{s,-\frac{1}{2}}} 
\lesssim \parallel u_1\parallel_{{W}_p^{s,\frac{1}{2}}}\parallel u_2\parallel_{{W}_p^{s,\frac{1}{2}}} ,
 \end{equation}
i.e. we first prove Prop \ref{Bilinearlp} with no gain of $T$.
Then, it suffices to show
\begin{equation}\label{Bilinearlp2}
\parallel H(u_1,u_2)(n,\tau)\parallel_{{\hat{W}}_p^{s,-\frac{1}{2}}} \lesssim \|u_1\|_{b^{0}_{p,\infty}L^p_\tau}\|u_2\|_{b^{0}_{p,\infty}L^p_\tau},
 \end{equation}
 where $H(.,.)$ is defined by
$$
H(u_1,u_2)(n,\tau)=\frac{1}{2\pi}\sum_{n_1+n_2=n}\frac{|n|<n>^s}{<n_1>^s<n_2>^s}
\int_{\tau_1+\tau_2=\tau}\frac{f(n_1,\tau_1)g(n_2,\tau_2)d\tau_1}{\prod_{j=1}^{2}w(n_j,\tau_j)<\tau_j-n_j^3+\frac{1}{n_j}>^{\frac{1}{2}}}.
$$

Let $ \sigma(\tau,\tau_1,n,n_1) = \text{max}\{\mid \tau + m(n)\mid, \mid \tau_1 + m(n_1)\mid, \mid \tau - \tau_1 +  m(n-n_1)\mid\}$, then by Lemma \ref{lemma1}
 , we have
$ \sigma \gtrsim \mid nn_1n_2\mid$.\\
Part 1. First, we consider the $\hat{X}_p^{0,-\frac{1}{2}}$ part of the $\hat{W}_p^{0,-\frac{1}{2}}$\\
Case (1): $\sigma = \mid \tau + m(n)\mid$. Suppose that $\mid n_1\mid \geq \mid n_2 \mid$. For fixed $n \neq 0$ and $ \tau$, let 
$\lambda = \frac{\tau-n^3+\frac{1}{n}}{3n}$ and define
$
B_{n,\tau} = \{n_1 \in \mathbb{Z}: \mid n_1 - r_j\mid \geq 1, j=1,2 r_j ~\text{is a real root
 of} ~p_{n,\lambda}(n_1):=\lambda + n_1(n-n_1),
 or ~ r_j = \frac{1}{2}n ~\text{if no real root} \}
$\\
On $B_{n,\tau}^{c}$ we have 
$$
\frac{|n|<n>^s}{<n_1>^s<n_2>^s\sigma^{\frac{1}{2}}}\lesssim \frac{1}{\mid n_2\mid^{\frac{1}{2}+s}}.
$$
By lemma \ref{int} we have
Note that for fixed $n$ and $\tau$ there are at most four values of $n_1 \in  B_{n,\tau}^{c}$,  i.e. the summation
over $n_1$ can be replace by the $L^p _{n_1}$ norm. Then, by H$\ddot{o}$lder inequality, we have
\begin{align}
 \text{LHS of} ~(\ref{Bilinearlp2}) &\lesssim \sup_{j}\|\sum_{n_1+n_2=n}\frac{w(n,\tau)}{<n_2>^{\frac{1}{2}+s}}
\int_{\tau_1+\tau_2=\tau}\frac{f(n_1,\tau_1)g(n_2,\tau_2)d\tau_1}{\prod_{j=1}^{2}<\tau_j-n_j^3+\frac{1}{n_j}>^{\frac{1}{2}}}\|_{L^p_{|n|\sim2^j}L^p_\tau}\nonumber\\
&\lesssim \sup_{j}\| \frac{w(n,\tau)}{<n_2>^{\frac{1}{2}+s}}\|f(n_1,)\|_{L^p_\tau}\|g(n_2,)\|_{L^p_\tau}\|_{L^p_{|n|\sim 2^j}L^p_{n_1}}.\nonumber
\end{align}
Note that $ w(n,\tau) \lesssim <n_2>^\delta$ since $|n_1|\geq|n_2|$.\\
If $|n_1|\gg|n_2|$ and $|n| \sim 2^j$, then we have $ |n_1| \sim 2^k$ where $\mid k-j\mid \leq 5$.
\begin{align}
 \text{LHS of} ~(\ref{Bilinearlp2})&\lesssim \sup
\bigg(\sum_{\mid k-j\mid \leq 5}\sum_{|n_1| \sim 2^k}\sum_{l}\sum_{\mid n_2\mid \sim 2^l}<n_2>^{(-\frac{1}{2}-s+\delta)p}
\|f(n_1,)\|^p_{L^p_\tau}\|g(n_2,)\|^p_{L^p_\tau}\bigg)^{\frac{1}{p}}\nonumber\\
&\lesssim \sum_{l}2^{(-\frac{1}{2}-s+\delta)pl}\sup_{k}
\parallel f \parallel_{L^p_{|n|\sim 2^K}L^p_{\tau}}\sup_{l}\parallel g \parallel_{L^p_{|n|\sim 2^l}L^p_{\tau}}
\lesssim \|f\|_{b^{0}_{p,\infty}L^p_\tau}\|g\|_{b^{0}_{p,\infty}L^p_\tau}\nonumber
\end{align}
by taking $\delta >0$ suifficiently small that $ -\frac{1}{2}-s+\delta <0$. Similarly, if $\mid n_1\mid \sim \mid n_2\mid$ and   $\mid n_2\mid \sim 2^l$, 
then we have $\mid n_1\mid \sim 2^k$ where $\mid k-l\mid \leq 5$,
\begin{align}
 \text{LHS of}~(\ref{Bilinearlp2})&\lesssim \bigg(
\sum_{\mid k-l\mid \leq 5}\sum_{|n_1| \sim 2^k}\sum_{l}\sum_{\mid n_2\mid \sim 2^l}<n_2>^{(-\frac{1}{2}-s+\delta)p}
\|f(n_1,)\|^p_{L^p_\tau}\|g(n_2,)\|^p_{L^p_\tau}\bigg)^{\frac{1}{p}}\nonumber\\
&\lesssim \sum_{l}2^{(-\frac{1}{2}-s+\delta)pl}\sup_{k}
\parallel f \parallel_{L^p_{|n|\sim 2^K}L^p_{\tau}}\sup_{l}\parallel g \parallel_{L^p_{|n|\sim 2^l}L^p_{\tau}}
\lesssim \|f\|_{b^{0}_{p,\infty}L^p_\tau}\|g\|_{b^{0}_{p,\infty}L^p_\tau}.\nonumber
\end{align}
Now on $B_{n,\tau}$ we have $<\tau - m(n)> \gtrsim <n><\lambda + n_1(n-n_1)>$. Also recall that $ w(n,\tau) \lesssim <\tau-m(n)>^{0^+}$.\\
Moreover, $ <\tau-m(n)>^{0^+} \lesssim \max (<n>,<n_2>, <\tau-m(n)+3nn_1n_2>)^{0^+}$ then we have
\begin{equation}\label{w}
 w(n,\tau) \lesssim (<n_2><\tau-m(n)+3nn_1n_2>)^{0^+}.
\end{equation}
By applying H\"{o}lder inequality and proceeding as before, we obtain:\\
LHS of $(\ref{Bilinearlp2}) \lesssim M \sup_j\parallel <n_2>^{0-}\|f(n_1,)\|_{L^p_\tau}\|g(n_2,)\|_{L^p_\tau}\|_{L^p_{|n|\sim 2^j}L^p_{n_1}}
\lesssim M \|f\|_{b^{0}_{p,\infty}L^p_\tau}\|g\|_{b^{0}_{p,\infty}L^p_\tau},$
where 
$$
M= \sup_{n,\tau}\parallel \frac{w(n,\tau)}{<n_2>^{\frac{1}{2}+s-}<\tau-m(n)+3nn_1n_2>^{1-\frac{1}{p^{\prime}}}}\parallel_{L^{p^{\prime}}_{n_1}}.
$$
Is suffice now to show that $M <\infty$, by (\ref{w}) and Lemma \ref{sum}, we have 
$$
M^{p^{\prime}} \lesssim \sup_{n,\tau}\frac{1}{<n>^{p^{\prime}-1_{-}}}
\sum_{n_2}\frac{1}{<n_2>^{(\frac{1}{2}+s-)p^{\prime}}<\lambda + n_1(n-n_1)>^{p^{\prime}-1-}} < \infty,
$$
since $ (\frac{1}{2}+s-)p^{\prime} + 2(p^{\prime}-1)>1$ for $p = 2+<4$ ans $sp = -1-$.\\
Now, assume $ \sigma = <\tau_2 - m(n_2)>$ (By symmetry this takes care of the case when $\sigma = <\tau_1 - m(n_1)>$).
Note that we have $ w(n,\tau) \lesssim <\tau-m(n)>^{0^+} $. Thus for this case by duality it suffices to show that:
\begin{align}\label{Bilinearlp3}
\sum_{l}\parallel\sum_{n}\frac{|n|<n>^s}{<n_1>^s<n_2>^s}&\frac{1}{w(n_2,\tau_2)<\tau_2-m(n_2)>^{\frac{1}{2}}}\nonumber\\
&\times\int \frac{f(n_1,\tau_1)h(n,\tau)d\tau}{<\tau_1-n_1^3+\frac{1}{n_1}>^{\frac{1}{2}}<\tau-n^3+\frac{1}{n}>^{\frac{1}{2}-}}
\parallel_{L^{p^{\prime}}_{|n_2|\sim 2^l}L^p_{\tau_2}}\nonumber\\
&\lesssim \sup_{k}\|f\|_{L^p_{|n_1|\sim 2^k}L^p_{\tau_1}}\sum_{j}\|h\|_{L^{p^{\prime}}_{\tau}L^{p^{\prime}}_{|n|\sim 2^j}}.
 \end{align}
For fixed $n_2 \neq 0$ and $\tau_2$, let $\lambda = \frac{\tau_2-n_2^3+\frac{1}{n_2}}{3n_2}$ and define

$
B_{n_2,\tau_2} = \{n_1 \in \mathbb{Z}: \mid n - r_j\mid \geq 1, j=1,2 r_j ~\text{is a real root
 of} ~p_{n_2,\lambda}(n_1):=\lambda + n(n_2-n),
 or ~ r_j = \frac{1}{2}n_2 ~\text{if no real root} \}
$\\

Case (2): $\sigma = <\tau_2 + m(n_2)>$ and $ \mid n_1\mid \gtrsim \mid n_2\mid$. \\
Subcase (2.a):On $B_{n_2,\tau_2}^c$ \\
First, suppose $<\tau_2-n_2^3 + \frac{1}{n_2} -3nn_1n_2> \gtrsim <n_2>^{\frac{1}{100}}$. Thus by Lemma \ref{int}, we have :
\begin{equation}
 \parallel <\tau_1-m(n_1)>^{-\frac{1}{2}+\alpha}<\tau_1-m(n)>^{-\frac{1}{2}}\parallel_{L^p_\tau}
\lesssim <\tau_2-n_2^3 + \frac{1}{n_2} -3nn_1n_2> ^{-\frac{1}{2}+\alpha+}\lesssim <n_2>^{-\frac{1}{100}(\frac{1}{2}-\alpha)+}
\end{equation}
for $ \alpha > 0$. Then by H\"{o}lder inequality in $\tau$ followed by Young and H\"{o}lder inequalities, we have
\begin{align}
\rVert \int &\frac{f(n_1,\tau_1)h(n,\tau)d\tau}{<\tau_1-n_1^3+\frac{1}{n_1}>^{\frac{1}{2}}<\tau-n^3+\frac{1}{n}>^{\frac{1}{2}-}}
\parallel_{L^p_{_{\tau_2}}} \lesssim <n_2>^{-\frac{1}{100}(\frac{1}{2}-\alpha)+}\|\frac{f(n_1,\tau_1)}{<\tau_1-n_1^3+\frac{1}{n_1}>^{\alpha}}
h(n,\tau)\|_{L^{p^{\prime}}_{\tau_2,\tau}}\nonumber\\
&\leq <n_2>^{-\frac{1}{100}(\frac{1}{2}-\alpha)+}\|<\tau_1-n_1^3+\frac{1}{n_1}>^{-\alpha}\|_{L^{\frac{p}{p-2}}_{\tau_1}}
\|f(n_1,)\|_{L_{\tau_1}^p}\|h(n,)\|_{L_{\tau}^{p^{\prime}}}\nonumber
 \end{align}
for fixed $n$ and $n_1$. By choosing $ \alpha> \frac{p-2}{p} = 0+$, we obtain\\
 $\|<\tau_1-n_1^3+\frac{1}{n_1}>^{-\alpha}\|_{L^{\frac{p}{p-2}}_{\tau_1}} < C < \infty$, independently of $n_1$.\\
Note that if $|n| \sim 2^j$ and $\mid n_2\mid \sim 2^l$, theen we have $\mid n_1\mid \sim 2^k$ where $\mid k-j \mid \leq 5$ or $\mid k- l \mid\leq 5$ since
$n=n_1+n_2$ and $\mid n_1 \mid \geqslant \mid n_2\mid$. For fixed $n_2$ and $\tau_2$ there are at most four values of $n \in B_{n_2,\tau_2}^c$ i.e the summation
over $n$ can be replace by the $L_n^{p^{\prime}}$ norm.\\
By H\"{o}lder inequality in $n_2$ after switching the order of summations, 
\begin{align}\label{Bilinearlp4}
 \text{LHS of}~ (\ref{Bilinearlp2})&\lesssim\sum_{l}\left\|<n_2>^{-\frac{1}{2}-s-\frac{1}{100}(\frac{1}{2}-\alpha)+}\left\|f(n_1,)\right\|_{L^p_{\tau_1}}
 \left\|h(n,)\right\|_{L^{\prime{p}_\tau}}\right\|_{L^{p^{\prime}_{\left|n_2\right|\sim2^l}}L^{p^{\prime}_\tau}}\nonumber\\
 &\bigg(\sum_{l}(2^l)^{0-}\bigg)\sup_{l}\bigg(\sum_j\sum_{\left|n\right|\sim2^j}\left\|
 <n_2>^{-\frac{1}{2}-s-\frac{1}{100}(\frac{1}{2}-\alpha)+}\right\|^{p^{\prime}}_{L^{\frac{p}{p-2}}_{n_2}}\nonumber\\
&\times\left\|f(n-n_2,)\right\|^{p^{\prime}}_{L^p_{\left|n_2\right|\sim2^l}L^p_{\tau_1}}\left\|h(n,)\right\|^{p^{\prime}}_{L^{p^{\prime}}_\tau}\bigg)^{\frac{1}{p^{\prime}}}\nonumber\\
 &\lesssim K\left\|f\right\|_{b^0_{p,\infty}L^p_{\tau}}\left\|h\right\|_{b^0_{p^{\prime},1}L^{p^{\prime}}_{\tau}},
\end{align}
where $K = \left\| <n_2>^{-\frac{1}{2}-s-\frac{1}{100}(\frac{1}{2}-\alpha)+}\right\|_{L_{n_2}^{\frac{p}{p-2}}}< \infty$(for $p < \frac{2_{-}}{1+\frac{1}{100}+}$). Note that we did not make use of $w(n_2,\tau_2)$ in this case.\\
Now, suppose   $<\tau_2-n_2^3 + \frac{1}{n_2} -3nn_1n_2> \ll <n_2>^{\frac{1}{100}}$, then we have $w(n_2,\tau_2)\sim <n>^{\delta}$ since $|n_1|\geq |n_2|$ implies $|n|\lesssim |n_1|$.If $|n_2|\lesssim |n|^{100}$, then $w(n_2,\tau_2) \geq <n_2>^{\frac{\delta}{100}}$. Otherwise, we have $|n|^{100}\lesssim |n_1|$, then we have
$$
\frac{|n|<n>^s}{<n_1>^s<n_2>^s}\frac{1}{\sigma^{\frac{1}{2}}} \lesssim \frac{1}{<n_1>^{(\frac{1}{2}+s)\frac{99}{100}}<n_2>^{\frac{1}{2}+s}}\lesssim
\frac{1}{<n_2>^{\frac{1}{2}+s+\epsilon}},
$$
for some $\epsilon >0$. Hence, we obtain a small power of $<n_2>$ in either case.\\
Subcase (2.b): On $B_{n_2,\tau_2}$. In this case we have
$$
<n_2><\frac{\tau_2-n_2^3+\frac{1}{n_2}}{3n_2} +n(n_2-n)>\lesssim <\tau_2-n_2^3+\frac{1}{n_2} -3nn_1n_2>. 
$$
Then with Holder inequality, we have
\begin{align}
\int \frac{f(n_1,\tau_1)h(n,\tau)d\tau}{<\tau_1-n_1^3+\frac{1}{n_1}>^{\frac{1}{2}}<\tau-n^3+\frac{1}{n}>^{\frac{1}{2}-}}&\lesssim
 <n_2>^{-\frac{1}{2}+\alpha+}<\frac{\tau_2-n_2^3+\frac{1}{n_2}}{3n_2} +n(n_2-n)>^{-\frac{1}{2}+\alpha+}\nonumber\\
 &\times\left\|\frac{f(n_1,\tau_1)}{<\tau_1 - n_1^3 +\frac{1}{n_1}>^{\alpha}}h(n,\tau)\right\|_{L^{p^{\prime}}_{\tau}}\nonumber
\end{align}
for fixed $n,~n_2,$and $\tau_2$.\\
Note that in this case we have:
$$
\frac{|n|<n>^s}{<n_1>^s<n_2>^s}\frac{1}{\sigma^{\frac{1}{2}}}\lesssim\frac{1}{<n_2>^{\frac{1}{2}+s}}.
$$
Now by Holder inequality in $n$ and $\tau_1$, we obtain that:
\begin{align}
\text{LHS of}~&  (\ref{Bilinearlp2})\lesssim \sum_{l}(2^l)^{0_{-}}K_1
\left\| <n_2>^{-1+\alpha-s+}\left\|\frac{f(n_1,\tau_1)}{<\tau_1 - n_1^3 +\frac{1}{n_1}>^{\alpha}}h(n,\tau)\right\|_{L^{p^{\prime}}_{\tau_2,\tau}}\right\|_{L^{p^{\prime}}_{\left|n_2\right|\sim 2^l}L^{p^{\prime}}_n}\nonumber\\
&\lesssim \sup_{l}K_1\left\| <n_2>^{-1+\alpha-s+}\left\|<\tau_1 - n_1^3 +\frac{1}{n_1}>^{-\alpha}\right\|_{L^{\frac{p}{p-2}}_{\tau_1}}
\left\|f(n_1,)\right\|_{L^p_{\tau_1}}\left\|h(n,)\right\|_{L^{p^{\prime}}_{\tau}}\right\|_{L^{p^{\prime}}_{\left|n_2\right|\sim 2^l}L^{p^{\prime}}_n},\nonumber
\end{align}
where $K_1 = \sup_{n_2,\tau_2}\left\|<\frac{\tau_2-n_2^3+\frac{1}{n_2}}{3n_2} +n(n_2-n)>^{-\frac{1}{2}+\alpha+}\right\|_{L^p_n} <\infty$ from Lemma \ref{sum}.
 We also have $\left\|<\tau_1 - n_1^3 +\frac{1}{n_1}>^{-\alpha}\right\|_{L^{\frac{p}{p-2}}_{\tau_1}} <\infty$, independently of $n_1$ as before.\\
Note that if  $|n| \sim 2^j$ and $|n_2|\sim 2^l$then we have $ |n_1| \sim 2^k$ where $\mid k-j\mid \leq 5$ or $\mid k-l\mid \leq 5$ since $n=n_1+n_2$ and $|n_2|\lesssim |n_1|$.
 Then, by Holder inequality in $n_2$, we have
\begin{align}
\text{LHS of}~ (\ref{Bilinearlp2}) &\lesssim K_2 \sup_{l}\bigg(\sum_{j}\sum_{\left|n\right|\sim 2^j}
\left\|f(n_1,)\right\|^{p^{\prime}}_{L^p_{\left|n_1\right|\sim 2^k}}\left\|h(n,)\right\|^{p^{\prime}}_{L^{p^{\prime}}_{\tau}}\bigg)^{\frac{1}{p^{\prime}}}\nonumber\\
&\lesssim K_2 \left\|f\right\|_{b^0_{p,\infty}L^p_{\tau}}\left\|h\right\|_{b^0_{p^{\prime},1}L^{p^{\prime}_{\tau}}},\nonumber
\end{align}
where $K_2 = \left\|<n_2>^{-1+\alpha -s+}\right\|_{L^{\frac{p}{p-2}}_{n_2}} < \infty $ since $(1-\alpha+s)\frac{p}{p-2} >1$.\\
Case(3):$\sigma = <\tau_2 - n_2^3 + \frac{1}{n_2}> $ and $\mid n_1\mid \ll \mid n_2\mid \Longrightarrow \mid n_1\mid  \ll \mid n_2\mid  \sim \mid n\mid$.\\
In this case, we have:
\begin{equation}\label{Bilinearlp5}
\frac{|n|<n>^s}{<n_1>^s<n_2>^s}\frac{1}{\sigma^{\frac{1}{2}}}\lesssim\frac{1}{<n_1>^{\frac{1}{2}+s}}.
\end{equation}
Subcase (3.a): On $B^c_{n_2,\tau_2}$.\\
If $<\tau_2-n_2^3 + \frac{1}{n_2} -3nn_1n_2> \gtrsim <n_2>^{\frac{1}{100}}$, 
then we have $ <\tau_2-n_2^3 + \frac{1}{n_2} -3nn_1n_2>  \gg <n_1>^{\frac{1}{100}}$.
As in Subcase (2.a) we have:
\begin{equation}
 \parallel <\tau_1-m(n_1)>^{-\frac{1}{2}+\alpha}<\tau_1-m(n)>^{-\frac{1}{2}}\parallel_{L^p_\tau}
\lesssim <\tau_2-n_2^3 + \frac{1}{n_2} -3nn_1n_2> ^{-\frac{1}{2}+\alpha+}\lesssim <n_1>^{-\frac{1}{100}(\frac{1}{2}-\alpha)+}
\end{equation}
and 

\begin{align}
\rVert \int &\frac{f(n_1,\tau_1)h(n,\tau)d\tau}{<\tau_1-n_1^3+\frac{1}{n_1}>^{\frac{1}{2}}<\tau-n^3+\frac{1}{n}>^{\frac{1}{2}-}}
\parallel_{L^p_{_{\tau_2}}} \lesssim <n_1>^{-\frac{1}{100}(\frac{1}{2}-\alpha)+}\|\frac{f(n_1,\tau_1)}{<\tau_1-n_1^3+\frac{1}{n_1}>^{\alpha}}
h(n,\tau)\|_{L^{p^{\prime}}_{\tau_2,\tau}}\nonumber\\
&\leq <n_1>^{-\frac{1}{100}(\frac{1}{2}-\alpha)+}\|<\tau_1-n_1^3+\frac{1}{n_1}>^{-\alpha}\|_{L^{\frac{p}{p-2}}_{\tau_1}}
\|f(n_1,)\|_{L_{\tau_1}^p}\|h(n,)\|_{L_{\tau}^{p^{\prime}}}.\nonumber
 \end{align}
Note that if $|n_2| \sim 2^l$ then we have $\mid n_1\mid \sim 2^k$ and  $\mid n\mid \sim 2^j$ where $k=0,...,l$ $\mid j- l \mid\leq 5$, then by Holder inequality

\begin{align}\label{Bilinearlp4}
 \text{LHS of}~ (\ref{Bilinearlp2})&\lesssim\sum_{k}\left\|<n_1>^{-\frac{1}{2}-s-\frac{1}{100}(\frac{1}{2}-\alpha)+}\left\|f(n_1,)\right\|_{L^p_{\tau_1}}
 \left\|h(n,)\right\|_{L^{\prime{p}_\tau}}\right\|_{L^{p^{\prime}_{\left|n_1\right|\sim2^k}}L^{p^{\prime}_\tau}}\nonumber\\
 &\bigg(\sum_{k}(2^k)^{0-}\bigg)\sup_{l}\bigg(\sum_j\sum_{\left|n\right|\sim2^j}\left\|
 <n_1>^{-\frac{1}{2}-s-\frac{1}{100}(\frac{1}{2}-\alpha)+}\right\|^{p^{\prime}}_{L^{\frac{p}{p-2}}_{n_1}}\nonumber\\
&\times\left\|f(n-n_2,)\right\|^{p^{\prime}}_{L^p_{\left|n_1\right|\sim2^k}L^p_{\tau_1}}\left\|h(n,)\right\|^{p^{\prime}}_{L^{p^{\prime}}_\tau}\bigg)^{\frac{1}{p^{\prime}}}\nonumber\\
 &\lesssim K_3\left\|f\right\|_{b^0_{p,\infty}L^p_{\tau}}\left\|h\right\|_{b^0_{p^{\prime},1}L^{p^{\prime}}_{\tau}},
\end{align}
where $K_3 = \left\| <n_1>^{-\frac{1}{2}-s-\frac{1}{100}(\frac{1}{2}-\alpha)+}\right\|_{L_{n_1}^{\frac{p}{p-2}}}< \infty$(for $p < \frac{2_{-}}{1+\frac{1}{100}+}$).\\
Now, suppose   $<\tau_2-n_2^3 + \frac{1}{n_2} -3nn_1n_2> \ll <n_2>^{\frac{1}{100}}$,
 then we have $w(n_2,\tau_2)\sim <n_1>^{\delta}$ since $|n_1|\ll |n|$. This extra gain of $<n_1>^{\delta}$ in the denominateur of (\ref{Bilinearlp2}) 
lets us proceed as before.\\
Subcase (3.b): On $B_{n_2,\tau_2}$. In this case 

$$
<n_2><\frac{\tau_2-n_2^3+\frac{1}{n_2}}{3n_2} +n(n_2-n)>\lesssim <\tau_2-n_2^3+\frac{1}{n_2} -3nn_1n_2>. 
$$ 
 Now, we have 

\begin{align}
\int \frac{f(n_1,\tau_1)h(n,\tau)d\tau}{<\tau_1-n_1^3+\frac{1}{n_1}>^{\frac{1}{2}}<\tau-n^3+\frac{1}{n}>^{\frac{1}{2}-}}&\lesssim
 <n_2>^{-\frac{1}{2}+\alpha+}<\frac{\tau_2-n_2^3+\frac{1}{n_2}}{3n_2} +n(n_2-n)>^{-\frac{1}{2}+\alpha+}\nonumber\\
 &\times\left\|\frac{f(n_1,\tau_1)}{<\tau_1 - n_1^3 +\frac{1}{n_1}>^{\alpha}}h(n,\tau)\right\|_{L^{p^{\prime}}_{\tau}}\nonumber
\end{align}
then using (\ref{Bilinearlp5}) we obtain

\begin{align}
\text{LHS of}~ (\ref{Bilinearlp2}) &\lesssim K_2 \sup_{l}\bigg(\sum_{j}\sum_{\left|n\right|\sim 2^j}
\left\|f(n_1,)\right\|^{p^{\prime}}_{L^p_{\left|n_1\right|\sim 2^k}}\left\|h(n,)\right\|^{p^{\prime}}_{L^{p^{\prime}}_{\tau}}\bigg)^{\frac{1}{p^{\prime}}}\nonumber\\
&\lesssim K_2 \left\|f\right\|_{b^0_{p,\infty}L^p_{\tau}}\left\|h\right\|_{b^0_{p^{\prime},1}L^{p^{\prime}_{\tau}}}.\nonumber
\end{align}
Part 2. Now we consider the $\hat{Y}_p^{0,-1}$ part of the $\hat{W}_p^{0,-\frac{1}{2}}$ norm. Define the bilinear operator $H_{\theta,b}(,)$ by
\begin{align}
H_{\theta,b}(f,g)(n,\tau)=\frac{1}{2\pi}\sum_{n_1+n_2=n}&\frac{1}{<\tau-m(n)>^{\theta}}\frac{|n|<n>^s}{<n_1>^s<n_2>^s}\nonumber\\
&\times\int_{\tau_1+\tau_2=\tau}\frac{f(n_1,\tau_1)g(n_2,\tau_2)d\tau_1}{\prod_{j=1}^{2}w(n_j,\tau_j)<\tau_j-m(n_j)>^{b}}.\nonumber
\end{align}
If $ \sigma = <\tau_1-m(n_1)>$ or $<\tau_2-m(n_2)>$, then by Holder inequality, we have
\begin{align}
 \text{LHS of}~ (\ref{Bilinearlp}) &= \sup_{j}\rVert H_{-1,\frac{1}{2}}(f,g)(n,\tau)\Arrowvert_{L^p_{|n| \sim 2^j}L^1_\tau}\nonumber\\
&\leq \sup_{j}\| \Arrowvert <\tau-m(n)>^{-\frac{1}{2} - \epsilon}
\Arrowvert_{L^{p^{\prime}}_\tau}\Arrowvert H_{-\frac{1}{2}+\epsilon,\frac{1}{2}}(f,g)(n,\tau)\Arrowvert_{L^p_\tau}\rVert_{L^p_{|n|\sim 2^j}}\nonumber\\
&\lesssim \sup_{j}\|H_{-\frac{1}{2}+\epsilon,\frac{1}{2}}(f,g)(n,\tau)\|_{L^p_{|n|\sim 2^j}L^p_\tau},\nonumber
\end{align}
where we choose $\epsilon >0$ such that $(\frac{1}{2}+\epsilon)p^{\prime}>1$. Then the proof reduces to Cases (2) and (3), where $<\tau -m(n)>^{\frac{1}{2}}$
is replaced by $<\tau -m(n)>^{\frac{1}{2}-\epsilon}$. Note that this does not affect the argument in Cases (2) and (3).\\
Now, assume $\sigma = <\tau-m(n)>$. If $\max(<\tau_1-m(n_1)>,<\tau_2-m(n_2)>)\gtrsim <\tau-m(n)>^{\frac{1}{100}}$, then By Holder inequality, we have
\begin{align}
 \text{LHS of}~ (\ref{Bilinearlp}) &\leq \sup_{j}\| \Arrowvert <\tau-m(n)>^{-\frac{1}{2} - \epsilon}\parallel_{L^{p^{\prime}}_\tau}
\Arrowvert H_{-\frac{1}{2},\frac{1}{2}-100\epsilon}(f,g)(n,\tau)\Arrowvert_{L^p_\tau}\rVert_{L^p_{|n|\sim 2^j}}\nonumber\\
&\lesssim\sup_{j}\Arrowvert H_{-\frac{1}{2},\frac{1}{2}-100\epsilon}(f,g)(n,\tau)\Arrowvert_{{L^p_\tau}L^p_{|n|\sim 2^j}}.\nonumber
\end{align}
Then, the proof reduces to Case (1) with $<\tau_j -m(n_j)>^{\frac{1}{2}}$ replaced by $<\tau_j -m(n_j)>^{\frac{1}{2}-100\epsilon}$,
which does not affect the argument.\\
Lastly, if $\max(<\tau_1 -m(n_1)>,<\tau_2-m(n_2)>) \ll <\tau-m(n)>^{\frac{1}{100}}$, then by Holder inequality, we have
\begin{align}
 \text{LHS of}~ (\ref{Bilinearlp}) &\leq \sup_{j}\parallel \parallel<\tau-m(n)>^{-\frac{1}{2}}1_{K(n)}(\tau-m(n))\parallel_{L^{p^{\prime}}}
\Arrowvert H_{-\frac{1}{2},\frac{1}{2}}(f,g)(n,\tau)\Arrowvert_{{L^p_\tau}}\parallel_{L^p_{|n|\sim 2^j}}\nonumber\\
&\lesssim \sup_{j}\Arrowvert H_{-\frac{1}{2},\frac{1}{2}}(f,g)(n,\tau)\Arrowvert_{{L^p_\tau}L^p_{|n|\sim 2^j}},\nonumber
\end{align}
where the second inequality follows from Lemma \ref{intOmega} since $-\frac{1}{2}p^{\prime} = -1+$. Once again, the proof reduces to Case (1).\\

Part 3: Now we discuss how to gain a small power of $T$ in (\ref{Bilinearlp}). From the two first parts we have showed:
\begin{equation}\label{Bilinearlp6}
\parallel \partial_x(u_1u_2)\parallel_{{W}_p^{s,-\frac{1}{2}}} 
\lesssim \parallel \hat{u_1}\parallel_{{\hat{X}}_p^{s,b}}\parallel w\hat{ u_1}\parallel_{{\hat{X}}_p^{s,\frac{1}{2}}}+
\parallel w\hat{u_2}\parallel_{{\hat{X}}_p^{s,\frac{1}{2}}}\parallel\hat{ u_2}\parallel_{{\hat{X}}_p^{s,b}}
 \end{equation}
for some $b<\frac{1}{2}$.
Thus, (\ref{Bilinearlp}) follows once we prove:
\begin{equation}\label{Bilinearlp7}
 \parallel \eta_{2T}(u)\parallel_{{X}_p^{s,b}} \lesssim T^{\theta} \parallel u \parallel_{{X}_p^{s,\frac{1}{2}}}. 
\end{equation}
By interpolation, we have
\begin{equation}\label{Bilinearlp7bis}
\parallel u \parallel_{{X}_p^{s,b}} \lesssim \parallel u \parallel^{\alpha}_{{X}_p^{s,0}}\parallel u \parallel^{1-\alpha}_{{X}_p^{s,\frac{1}{2}}}
\end{equation}
where 
$\alpha = 1-2b$.\\ Recall $\hat{\eta_{2T}}(\tau)=2T\hat{\eta}(2T\tau)$,
 then $\parallel \hat{\eta}_{2T} \parallel_{L^q_\tau}\sim T^{\frac{q-1}{q}}
\parallel\hat{\eta}\parallel_{L^q_\tau}\sim T^{\frac{q-1}{q}}$.\\
For fixed $n$, by Young and Holder inequalities, we have
$$
\parallel\hat{\eta_{2T}}*\hat{u}(n,)\parallel_{L^p_\tau} \leq \parallel\hat{\eta_{2T}}\parallel_{L^{p^{\prime}}_\tau}
\parallel\hat{u}(n,)\parallel_{L^{\frac{p}{2}_\tau}} \lesssim T^{\frac{p^{\prime}-1}{p^{\prime}}}\parallel<\tau-m(n)>^{-\frac{1}{2}}\parallel_{L^p_\tau}
\parallel<\tau-m(n)>^{\frac{1}{2}}\hat{u}(n,)\parallel_{L^p_\tau}.
$$
Hence for $p>2$, we have
\begin{equation}\label{Bilinearlp8}
\parallel u \parallel_{{X}_p^{s,0}}\lesssim T^{\frac{1}{p}}\parallel u \parallel_{{X}_p^{s,\frac{1}{2}}}.
\end{equation}
Then, (\ref{Bilinearlp7}) follows from (\ref{Bilinearlp7bis}) and (\ref{Bilinearlp8}). This completes the proof.

\end{document}